\documentclass{amsart}
\usepackage{amscd}
\usepackage{amssymb}
\usepackage{graphics}
\usepackage{pst-tree}
\usepackage[all]{xy}

\newtheorem{Thm}{Theorem}
\newtheorem{Prop}[Thm]{Proposition}
\newtheorem{Lem}[Thm]{Lemma}

\newtheorem*{Ack}{Acknowledgment}

\theoremstyle{definition}
\newtheorem{Def}[Thm]{Definition}
\newtheorem{Exa}[Thm]{Example}

\newcommand{\lora}{{\longrightarrow}}

\newcommand\de{\partial}

\newcommand{\Z}{\mathbb{Z}}

\newcommand{\End}{\mathcal E\!nd}

\newcommand{\Det}{\mathrm{Det}}

\newcommand{\f}{
 \begin{pspicture}(0,0.1)(0.2,0.4)
 \psline[linewidth=1pt](0.1,0)(0.1,0.4)
 \end{pspicture}}
\newcommand{\da}{
 \begin{pspicture}(0,0.1)(0.2,0.4)
 \psline[linewidth=1pt, linestyle=dashed,dash=4pt 3pt](0.1,0)(0.1,0.4)
 \end{pspicture}}
\newcommand{\ff}{
 \begin{pspicture}(0,0)(0.2,0.2)
 \psline[linewidth=1pt](0.1,0)(0.1,0.2)
 \end{pspicture}}
\newcommand{\dda}{
 \begin{pspicture}(0,0)(0.2,0.2)
 \psline[linewidth=1pt, linestyle=dashed,dash=3pt 2pt](0.1,0)(0.1,0.2)
 \end{pspicture}}

\def\cmp#1,{{ Commun.\ Math.\ Phys.\ \bf #1},}
\def\jmp#1,{{ J.\ Math.\ Phys.\ \bf #1},}
\def\pl#1,{{ Phys.\ Lett.\ \bf #1},}
\def\npb#1,{{ Nucl.\ Phys.\ {\bf B #1}},}
\def\mpl#1,{{ Mod.\ Phys.\ Lett.\ \bf #1},}
\def\pr#1,{{ Phys.\ Rev.\ \bf #1},}
\def\prl#1,{{ Phys.\ Rev.\ Lett.\ \bf #1},}
\def\lmp#1,{{ Lett.\ Math.\ Phys.\ \bf #1},}
\def\jktr#1,{{ Jour.\ of Knot Theory and its Ramification \bf #1},}
\def\bams#1,{{ Bull.\ Amer.\ Math.\ Soc.\ \bf #1},}
\def\ja#1,{{ J.\ Algebra \bf #1},}
\def\ijm#1,{{ Illinois Journal of Mathematics \bf #1},}
\def\dmj#1,{{ Duke Math.\ J.\ \bf #1},}
\def\lnm#1,{{ Lecture Notes in Mathematics \bf #1},}
\def\mrl#1,{{ Math.\ Res.\ Lett.\ \bf #1},}

\begin{document}

\title{Homotopy Inner Products for Cyclic Operads}

\author[R.~Longoni]{Riccardo~Longoni}
\address{Dipartimento di Matematica ``G. Castelnuovo'',
Universit\`a di Roma ``La Sapienza'', Piazzale Aldo Moro, 2
I-00185 Roma, Italy}
\email{}

\author[T.~Tradler]{Thomas~Tradler}
\address{College of Technology of the City University
of New York, Department of Mathematics, 300 Jay Street, Brooklyn,
NY 11201, USA}
\email{ttradler@citytech.cuny.edu}

\begin{abstract}
We introduce the notion of homotopy inner products for any cyclic
quadratic Koszul operad $\mathcal O$, generalizing the construction
already known for the associative operad. This is done by defining a
colored operad $\widehat{\mathcal O}$, which describes modules over
$\mathcal O$ with invariant inner products. We show that
$\widehat{\mathcal O}$ satisfies Koszulness and identify algebras
over a resolution of $\widehat{\mathcal O}$ in terms of derivations
and module maps. As an application we construct a homotopy inner product over the commutative operad on the cochains of any Poincar\'e duality space.
\end{abstract}


\maketitle

\tableofcontents

\section{Introduction}

In \cite{GeK}, the notion of cyclic operads and invariant inner product for such operads was defined. A homotopy version of these inner products for the associative operad was given in \cite{Tr} and the starting point for a similar version for the commutative operad was considered in \cite{Ginot}. It is natural to ask for a generalization of these constructions applicable to any cyclic operad. This is what is done in this paper.

Starting with a cyclic operad $\mathcal O$, we use the notion of colored operads to incorporate the cyclic structure of $\mathcal O$ into the colored operad $\widehat{\mathcal O}$. Algebras over the colored operad $\widehat{\mathcal O}$ consist of pairs $(A,M)$, where $A$ is an algebra over $\mathcal O$ and $M$ is an $\mathcal O$-module over $A$ which has an invariant inner product. Section \ref{cyclic-op} is devoted to explicitly defining $\widehat{\mathcal O}$, and in the case that $\mathcal O$ is given by quadratic generators and relations, we give a description of $\widehat{\mathcal O}$ in terms of generators and relations coming from those of $\mathcal O$.

A major tool in the theory of operads is the notion of Koszul duality. Let us recall, that a (colored) operad $\mathcal P$ is called Koszul if there is a quasi-isomorphism of operads $\mathbf{D}(\mathcal P^!)\to \mathcal P$, where $\mathbf{D}(\mathcal P^!)$ denotes the dual operad (in the sense of \cite[(3.2.12)]{GK}) on the dual quadratic operad $\mathcal P^!$. This implies that the notion of algebras of $\mathcal P$ has a canonical infinity version given by algebras over $\mathcal P_\infty:=\mathbf{D}(\mathcal P^!)$. Our main theorem states, that the Koszulness property is preserved when going from $\mathcal O$ to  $\widehat{\mathcal O}$.
\setcounter{Thm}{8} \begin{Thm}
Let $\mathcal O$ be a cyclic quadratic operad. If $\mathcal O$ is Koszul, then so is $\widehat{\mathcal O}$, i.e. we have a resolution $\widehat{\mathcal O}_\infty:=\mathbf{D}(\widehat{\mathcal O^!})$ of $\widehat{\mathcal O}$.
\end{Thm} \setcounter{Thm}{0}
The proof of this theorem will be given in section \ref{quadrat-koszul}. Theorem \ref{O_hat_Koszul} justifies the concept of the infinity version of algebras and modules with invariant inner products over cyclic operads $\mathcal O$ as algebras over the operad $\widehat{\mathcal O}_\infty=\mathbf{D}(\widehat{\mathcal O^!})$. The concept of algebras over the operad $\widehat{\mathcal O}_\infty$ will be investigated in more detail in section \ref{homotop-ip}. In particular, we explicitly reinterpret in proposition \ref{O_hat_algebras} algebras over $\widehat{\mathcal O}_\infty$ in terms of derivations and module maps over free $\mathcal O^!$-algebras and modules.

Recall that the associative operad $\mathcal Assoc$, the commutative operad $\mathcal Comm$ and the Lie operad $\mathcal Lie$ are all cyclic quadratic Koszul operads, so that theorem \ref{O_hat_Koszul} may be applied to all these cases. As a particularly application of infinity inner products, we consider the examples of the associative operad $\mathcal Assoc$ and the commutative operad $\mathcal Comm$, which have an interesting application to the chain level of a Poincar\'e duality space $X$. In \cite{TZ}, it was shown that the simplicial cochains $C^\ast(X)$ on $X$ with rational coefficients form an algebra over the operad $\widehat{\mathcal Assoc}_\infty$. This structure was then used in \cite{Tr2} and \cite{TZ2} to obtain string topology operations on the Hochschild cohomology, respectively the Hochschild cochain complex, of $C^\ast(X)$. Since the method of constructing the $\widehat{\mathcal Assoc}_\infty$ algebra on $C^\ast (X)$ easily transfers to the commutative case, we will show in section \ref{Comm-section} that $C^\ast(X)$ also forms a $\widehat{\mathcal Comm}_\infty$-algebra. We expect that this stronger algebraic structure should induce even more string topology operations, which take into account the commutative nature of the cochains of the space $X$. A first step in this direction was done in \cite{TZ3}, where the string topology operations for algebras over $\widehat{\mathcal Assoc}$ and $\widehat{\mathcal Comm}$ were investigated.

\begin{Ack}
We are grateful to Dennis Sullivan for many valuable suggestions and illuminating discussions. We also thank Domenico Fiorenza, Martin Markl, Jim Stasheff and Scott Wilson for useful comments and remarks regarding this topic. The second author was partially supported by the Max-Planck Institute in Bonn.
\end{Ack}

\section{$\widehat{\mathcal Comm}_\infty$ structure for Poincar\'e duality spaces}\label{Comm-section}

Before going into the details of the construction of homotopy inner products over a general cyclic quadratic operad $\mathcal O$, we give an application for the case of the commutative operad $\mathcal Comm$. More precisely, we show how a homotopy $\mathcal Comm$-inner product arises on the chain level of a Poincar\'e duality space $X$. In fact, the construction for the homotopy $\mathcal Comm$-algebra is taken from R. Lawrence and D. Sullivan's paper \cite{S} on the construction of local infinity structures. In \cite{TZ}, M. Zeinalian and the second author construct homotopy $\mathcal Assoc$-inner products on a Poincar\'e duality space $X$. The same reasoning may in fact be used to construct homotopy $\mathcal Comm$-inner products on $X$. The proof of the next proposition will be a sketch using these arguments.

\begin{Prop}\label{prop:comm-pd}
Let $X$ be a closed, finitely triangulated Poincar\'e duality space, such that the closure of every simplex is contractible. Denote by $C=C_\ast(X)$ the simplicial chains on $X$. Then its dual space $A:=C^*=Hom(C_*(X),k)$ has the structure of a $\widehat{\mathcal Comm}_\infty$ algebra, such that the lowest
multiplication is the symmetrized Alexander-Whitney multiplication
and the lowest inner product is given by capping with the fundamental cycle $\mu\in C$.
\end{Prop}
\begin{proof}[Proof (sketch)]
Let $\mathcal Lie$ denote the Lie-operad, $F_{\mathcal Lie}V=\bigoplus_{n\geq 1} (\mathcal Lie(n)\otimes V^{\otimes n})_{S_n}$ denote the free Lie algebra generated by $V$, and $F_{\mathcal Lie,V}W=\bigoplus_{n\geq 1} (\bigoplus_{k+l=n-1}\mathcal Lie(n)\otimes V^{\otimes k}\otimes W\otimes V^{\otimes l})_{S_n}$ the canonical module over $F_{\mathcal Lie}V$.
We will see in proposition \ref{O_hat_algebras} and example \ref{exa-comm}, that the required data for a homotopy $\mathcal Comm$-inner product consists of,
\begin{itemize}
\item
a derivation $d\in \mathrm{Der}(F _{\mathcal Lie}\,C[1])$ of degree
$1$, with $d^2=0$,
\item
a derivation $g\in \mathrm{Der}_d (F_{\mathcal Lie,\,C[1]}C[1])$ over
$d$ of degree $1$, with $g^2=0$, which imduces a derivation $h\in \mathrm{Der}_d (F_{\mathcal Lie,\,C[1]}C^*[1])$ over $d$ with $h^2=0$,
\item
a module map $f\in \mathrm{Mod}(F_{\mathcal Lie, C[1]}C^*, F_{\mathcal
Lie,C[1]}C[1])$ of degree $0$ such that $f\circ h = g \circ f$.
\end{itemize}

In order to construct the derivation $d\in \mathrm{Der}(F _{\mathcal Lie}\,C[1])$ with $d^2=0$, 
let $F_{\mathcal Lie}C[1]=L_1\oplus L_2\oplus\dots$, where $L_n=(\mathcal Lie(n)\otimes C[1]^{\otimes n})_{S_n}$, be the decomposition of $F_{\mathcal Lie}C[1]$
by the monomial degree in $C[1]$. Then, $d:F_{\mathcal Lie}C[1]\to F_{\mathcal Lie}C[1]$ is determined by maps $d=d_1+d_2+\dots$, where $d_i:C[1]\to L_i$ is lifted to $F_{\mathcal Lie}(C[1])$ as a derivation. Let $d_1$ be the differential on $C[1]$, and $d_2$ be the symmetrized Alexander-Whitney comultiplication. For the general $d_i$, we use the inductive
hypothesis that $d_1$, \dots, $d_{i-1}$ are local maps so that
$\nabla_i:=d_1+ \dots+d_{i-1}$ has a square $\nabla_i^2$ mapping only
into higher components $L_{i}\oplus L_{i+1}\oplus \dots$. Here, ``local'' means
that every simplex maps into the sub-Lie algebra of its closure. Now,
by the Jacobi-identity, it is true that $0=[\nabla_i, [\nabla_i,\nabla_i]]=
[d_1,e_i]\text{+ higher terms}$, where $e_i:C[1]\to L_i$ is the lowest
term of $[\nabla_i,\nabla_i]$. Thus $e_i$ is $[d_1,.]$-closed and
thus, using the contractibility hypothesis of the proposition, also locally
$[d_1,.]$-exact. These exact terms can be put together to give a map
$d_i$, so that $[d_1,d_i]$ vanishes on  $L_1\oplus \dots \oplus
L_{i-1}$ and equals $- 1/2\cdot e_i$ on $L_i$. In other words, $(d_1+
\dots+d_i)^2=1/2\cdot [d_1+ \dots+d_i,d_1+ \dots+d_i]=1/2\cdot [\nabla_i
,\nabla_i ]+ [d_1,d_i]+\text{higher terms}$, maps only into
$L_{i+1}\oplus L_{i+2}\oplus\dots$. This completes the inductive step,
and thus produces the wanted derivation $d$ on $F_{\mathcal Lie}(C[1])$. 

In a similar way, we may produce the derivation $g$ of $F_{\mathcal Lie,\,C[1]}C[1]$ over $d$, by  decomposing $F_{\mathcal Lie, C[1]} C[1]=L'_1\oplus L'_2\oplus\dots$, where $L'_n$ is given by the space $\left(\bigoplus_{k+l=n-1} \mathcal Lie(n) \otimes C[1]^{\otimes k} \otimes C[1] \otimes C[1]^{\otimes l}\right)_{S_{n}}$. With this notation, $g$ is written as a sum $g=g_1+g_2+\dots$, where $g_i:C[1]\to L'_i$ is lifted to $F_{\mathcal Lie, C[1]} C[1]$ as a derivation over $d$, and $(g_1+\dots+g_i)^2$ only maps into $L'_{i+1}\oplus L'_{i+2}\oplus\dots$.

Using a slight variation of the above method, we may also construct the wanted homotopy $\mathcal Comm$-inner product, i.e. the module map $f$ stated above. More precisely, we build a map $\chi:C[1]\to Mod(F_{\mathcal Lie, C[1]}C^*[1],F_{\mathcal Lie, C[1]}C[1])$, so that $\chi$ is a chain map under the differential $d_1$ on $C[1]$, and the differential $\delta(f)=f\circ h - (-1)^{|f|} g\circ f$ on $Mod(F_{\mathcal Lie, C[1]}C^*[1],F_{\mathcal Lie, C[1]}C[1])$. Since a module map is given by the components $M_n=\bigoplus_{k+l=n-2} \mathcal Lie(n)\otimes C[1]^{\otimes k}\otimes C[1]\otimes C[1]^{\otimes l}\otimes C[1]$, it is enough to construct $\chi$ as a sum $\chi=\chi_2+\chi_3+\dots$, where $\chi_i:C[1]\to M_i$. Now, the lowest component $\chi_2:C[1]\to C[1]\otimes C[1]$ is defined to be the symmetrized Alexander-Whitney comultiplication. For the induction, we assume that $\Upsilon_i:=\chi_2+\dots+\chi_{i-1}$ are local maps such that $D(\Upsilon_i):=\Upsilon_i\circ d_1-\delta\circ \Upsilon_i$ maps only into higher components $M_i\oplus M_{i+1}\oplus\dots$. Let $\epsilon_i:C[1]\to M_i$ be the lowest term of $D(\Upsilon_i)$. Since $D^2=0$ and $\delta$ has $d_1$ as its lowest component, we see that $\epsilon_i$ is $[d_1,.]$-closed, and by the hypothesis of the proposition also locally $[d_1,.]$-exact. These exact terms can be put together as before to produce a map $\chi_i$, so that now $\Upsilon_{i+1}:=\chi_2+\dots+\chi_i$ only maps into $M_{i+1}\oplus M_{i+2}\oplus\dots $. We therefore obtain the chain map $\chi$, and with this, we define the homotopy $\mathcal Comm$-inner product as $f:=\chi(\mu)\in  Mod(F_{\mathcal Lie, C[1]}C^*[1],F_{\mathcal Lie, C[1]}C[1]) $, where $\mu\in C$ denotes the fundamental cycle of the space $X$. Since $\mu$ is $d_1$-closed, it follows that $f\circ h-g\circ f=0$.
\end{proof}

\section{The operad $\widehat{\mathcal O}$}
\label{cyclic-op}

In this section, we define for any cyclic operad $\mathcal O$ the colored operad $\widehat{\mathcal O}$. In the case that $\mathcal O$ is cyclic quadratic, we give an explicit description of $\widehat{\mathcal O}$ in terms of generators and relations coming from generators and relations in $\mathcal O$.

We assume that the reader is familiar with the notion of operads, colored operads and cyclic operads.  For a good introduction to operads, we refer to \cite{Ad}, \cite{GK} and \cite{MSS}, for cyclic operads we recommend \cite{GeK} and \cite{MSS}. Colored operads were first introduced in 
\cite{BV} and appeared in many other places, see e.g. \cite{L} and \cite{BM}. Since in our case, we only need a special type of colored operad, it will be convenient to setup notation with the following definition.

As in \cite[(1.2.1)]{GK} and \cite[(1.1)]{GeK}, we assume throughout this paper that $k$ is a field of characteristic $0$. Note however, that for certain operads such as e.g. the associative operad, a more general setup is possible.
\begin{Def}\label{0/1-operad}
Let $\mathcal P$ be a 3-colored operad  in the category of (differential graded) vector spaces, where we use the three colors ``full", ``dashed" and ``empty", in symbols  written  $\f,\da,\varnothing$. This means that to each finite sequence of symbols $x_1, \dots, x_n, x\in \{\f,\da,\varnothing \}$, we have (differential graded) vector spaces $\mathcal P(x_1,\dots, x_n; x_0)$ over $k$, which label the operations with $n$ inputs with colors $x_1,\dots, x_n$, and one output with color $x$. The operad $\mathcal P$ comes with maps $\circ_i: \mathcal P(x_1,\dots,x_n;x) \otimes \mathcal P(y_1,\dots,y_m;x_i) \to \mathcal P(x_1,\dots,x_{i-1},y_1,\dots, y_m,x_{i+1},\dots, x_n;x)$ which label the composition in $\mathcal P$, and with an action of the symmetric group $S_n$, which, for $\sigma\in S_n$, maps $\mathcal P(x_1,\dots, x_n;x)\to \mathcal P(x_{\sigma(1)},\dots,x_{\sigma(n)};x)$. These maps have to satisfy the usual associativity and equivariance axioms of colored operads.

$\mathcal P$ is called a 0/1-operad if the color $\varnothing$ can appear
only as an output, and the only nontrivial spaces with one input are $\mathcal P(\f;\f)=k$ and
$\mathcal P(\da;\da)=k$. We assume furthermore, that there are fixed generators of the spaces $\mathcal P(\f;\f)=k$ and $\mathcal P(\da;\da)=k$.

Graphically, we represent $\mathcal P(x_1,\dots,x_n;x)$ by a tree with $n$ inputs and one output of the given color. Since the color $\varnothing$ cannot appear as an input, we may use the
following convention: we represent the output $\varnothing$ with a
blank line, i.e., with no line, and we say that the operation
``has no output''.
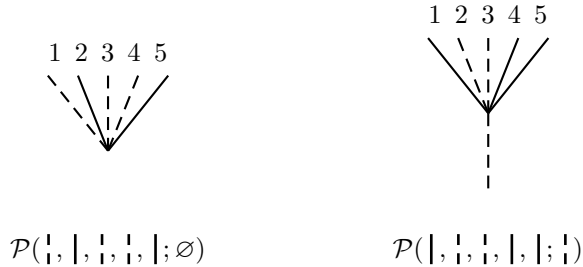
\[
\begin{pspicture}(0,.5)(4,4)
 \psline[linestyle=dashed, arrowsize=0.1, arrowinset=0](2,2)(1.2,3)
 \psline[arrowsize=0.1, arrowinset=0](2,2)(1.6,3)
 \psline[linestyle=dashed, arrowsize=0.1, arrowinset=0](2,2)(2,3)
 \psline[linestyle=dashed, arrowsize=0.1, arrowinset=0](2,2)(2.4,3)
 \psline[arrowsize=0.1, arrowinset=0](2,2)(2.8,3)
 \rput[b](2,.5){$\mathcal P (\da,\f,\da,\da,\f;\varnothing)$}
 \rput[b](2,3.2){$1\,\,\, 2\,\,\, 3\,\,\, 4\,\,\, 5$}
\end{pspicture}
\quad \quad \quad
\begin{pspicture}(0,0)(4,3.6)
 \psline[arrowsize=0.1, arrowinset=0](2,2)(1.2,3)
 \psline[linestyle=dashed, arrowsize=0.1, arrowinset=0](2,2)(1.6,3)
 \psline[linestyle=dashed, arrowsize=0.1, arrowinset=0](2,2)(2,3)
 \psline[arrowsize=0.1, arrowinset=0](2,2)(2.4,3)
 \psline[arrowsize=0.1, arrowinset=0](2,2)(2.8,3)
 \psline[linestyle=dashed, arrowsize=0.1, arrowinset=0](2,2)(2,1)
 \rput[b](2,0){$\mathcal P (\f,\da,\da,\f,\f;\da)$}
 \rput[b](2,3.2){$1\,\,\, 2\,\,\, 3\,\,\, 4\,\,\, 5$}
\end{pspicture}
\]
The canonical example of a 0/1-operad is the endomorphism 0/1-operad
given for $k$-vector spaces $A$ and $M$ by
\begin{align*}
\End^{A,M}(\vec X;\f)=& Hom(\text{tensor products of $A$ and $M$},A)\\
\End^{A,M}(\vec X;\da)=& Hom(\text{tensor products of $A$ and $M$},M)\\
\End^{A,M}(\vec X;\varnothing)=& Hom(\text{tensor products of $A$ and
$M$},k).
\end{align*}
With this notation $(A,M,k)$ is an algebra over the 0/1-operad $\mathcal P$ if there exists a 0/1-operad map  $\mathcal P \to \End^{A,M}$. By slight abuse of language we will also call the tuple $(A,M)$ an algebra over $\mathcal P$.
\end{Def}

It is our aim to define for each cyclic operad $\mathcal O$, the associated 0/1-operad $\widehat{ \mathcal O}$, which incorporates the cyclic structure as a colored operad. Before doing so, let us briefly recall the definition of a cyclic operad from \cite[Theorem (2.2)]{GeK}. 
\begin{Def}
Let $\mathcal O$ be a operad, i.e we have vector spaces $\mathcal O(n)$ for $n\geq 1$, composition maps $\circ_i:\mathcal O(n)\otimes \mathcal O(m)\to \mathcal O(n+m-1)$, and an $S_n$-action on $\mathcal O(n)$  for each $n$, satisfying the usual axioms, see \cite[(1.2.1)]{GK}. $\mathcal O$ is called {\em cyclic} if there is an action of
the symmetric group $S_{n+1}$ on $\mathcal O(n)$, which extends the
given $S_n$-action, and satisfies, for $1\in\mathcal O(1)$,
$\alpha\in \mathcal O(m)$, $\beta\in \mathcal O(n)$ the following
relations:
\begin{eqnarray}
\label{compos_cyclic1} \tau_2(1)&=&1,\\ \label{compos_cyclic2}
\tau_{m+n}(\alpha\circ_k \beta)&=&\tau_{m+1}(\alpha)\circ_{k+1}
\beta,\quad\quad\quad \text{ for } k<m \\ \label{compos_cyclic3}
\tau_{m+n}(\alpha\circ_m \beta)&=&\tau_{n+1}(\beta)\circ_1
\tau_{m+1}(\alpha),
\end{eqnarray}
where $\tau_{j}\in S_{j}$ denotes the cyclic rotation of $j$ elements
$\tau_{j} :=1\in\Z_{j}\subset S_{j}$.
\end{Def}

\begin{Def}
\label{def_O_hat}
Let $\mathcal O$ be a cyclic operad with $\mathcal O(1)=k$. For a sequence of $n$ input colors $\vec X=(x_1,\dots, x_n)$ and the output color $x$, where $x_1, \dots, x_n, x\in\{\f,\da,\varnothing\}$, let
\[
 \widehat{\mathcal O}(\vec X;x):=
   \begin{cases}
    \mathcal O(n) & \text{if } x \text{ is ``full'', and } \vec X=(\f,\ldots,\f),\\
    \mathcal O(n) & \text{if } x \text{ is ``dashed'', and $\vec X$ has
    exactly one ``dashed'' input}\\
    \mathcal O(n-1) & \text{if } x=\varnothing \text{ and
    $\vec X$ has exactly two ``dashed' inputs,} \\
    \{0\} & \text{otherwise}.
  \end{cases}
\]
The definition of $\widehat{\mathcal O}(\vec X,\varnothing)$ is
motivated by the idea that one considers trees with $n-1$ inputs and
one output, and then uses the $S_{n+1}$ action to turn this output into a new input:
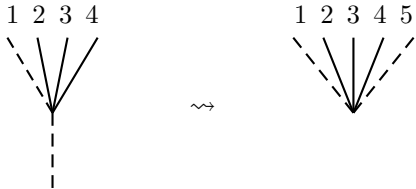
\[
\begin{pspicture}(0,0.8)(4,3.4)
 \psline[linestyle=dashed](2,2)(1.4,3)
 \psline(2,2)(1.8,3)
 \psline(2,2)(2.2,3)
 \psline(2,2)(2.6,3)
 \psline[linestyle=dashed](2,2)(2,1)
 \rput[b](2,3.2){$1\,\,\, 2\,\,\, 3\,\,\, 4$}
 \rput[b](4,2){$\rightsquigarrow$}
\end{pspicture}
\begin{pspicture}(0,0.8)(4,3.4)
 \psline[linestyle=dashed](2,2)(1.2,3)
 \psline(2,2)(1.6,3)
 \psline(2,2)(2,3)
 \psline(2,2)(2.4,3)
 \psline[linestyle=dashed](2,2)(2.8,3)
 \rput[b](2,3.2){$1\,\,\, 2\,\,\, 3\,\,\, 4\,\,\, 5$}
\end{pspicture}
\]
We define the $S_n$-action on $\widehat{\mathcal O}(\vec X;\f)$ and
$\widehat{\mathcal O}(\vec X;\da)$ as before by using the
$S_n$-action on $\mathcal O(n)$, and the $S_n$-action on
$\widehat{\mathcal O} (\vec X;\varnothing)$ by using the $S_n$-action
on $\mathcal O(n-1)$ given by the cyclicity of $\mathcal O$.

Diagrams with different positions of the two ``dashed'' inputs can be
mapped to each other using the action of the symmetric group. In
fact, as each $\sigma\in S_{n+1}$ induces an isomorphism which
preserves all the structure, any statement about diagrams with a
fixed choice of position of ``dashed'' inputs immediately carries
over to any other choice of positions of ``dashed'' inputs. We
therefore often restrict our attention to the choice where the two
``dashed'' inputs are at the far left and the far right, as shown in
the above picture.

It is left to define the composition. On
$\widehat{\mathcal O}(\vec X;\f)$ and $\widehat{\mathcal O}(\vec
X;\da)$, the composition is simply the composition in $\mathcal
O(n)$, so that it clearly satisfies associativity and equivariance.
If $|\vec X|=n+1$, then on $\widehat{\mathcal O}(\vec
X;\varnothing)=\mathcal O(n)$, the composition is predetermined on
the first $n$ components by the usual composition in $\mathcal O$. As
for the last component, we define
\begin{equation}\label{def_cyclic_compos}
\alpha\circ_{m+1} \beta:=\tau_{n+m} (\tau^{-1}_{m+1}(\alpha)\circ_m
\beta)\stackrel{\eqref{compos_cyclic3}}{=} \tau_{n+1}(\beta) \circ_1
\alpha
\end{equation}
\[
\begin{pspicture}(1,1.6)(10.5,5.6)
 \psline[linestyle=dashed](2,2)(1.2,2.9)
 \psline(2,2)(1.6,2.9)
 \psline(2,2)(2,2.9)
 \psline(2,2)(2.4,2.9)
 \psline[linestyle=dashed](2,2)(2.8,2.9)
 \rput[b](3,2){$\alpha$} \rput[b](2.4,3.4){$\beta$}
 \psline[linestyle=dashed](2.8,3)(2.8,3.5)
 \psline[linestyle=dashed](2.8,3.5)(3,4)
 \psline(2.8,3.5)(2.8,4)
 \psline(2.8,3.5)(2.6,4)
 \rput[b](3.5,3){$:=$}
 \psline[linestyle=dashed](5,2)(4.2,2.9)
 \psline[linestyle=dashed](4.2,3)(4.3,3.1)(5,3.25)(5.7,3.4)(5.8,3.5)
 \psline(5,2)(4.6,2.9)                   \psline(4.6,3)(4.2,3.5)
 \psline(5,2)(5,2.9)                     \psline(5,3)(4.6,3.5)
 \psline(5,2)(5.4,2.9)                   \psline(5.4,3)(5,3.5)
 \psline[linestyle=dashed](5,2)(5.8,2.9) \psline[linestyle=dashed](5.8,3)(5.4,3.5)
 \rput[b](6,2){$\alpha$} \rput[b](6.4,3.1){$\tau_{m+1}^{-1}$}
 \psline[linestyle=dashed](5.4,3.6)(5.4,4.1)
 \psline[linestyle=dashed](5.4,4.1)(5.6,4.5)
 \psline(5.4,4.1)(5.4,4.5)
 \psline(5.4,4.1)(5.2,4.5)
 \rput[b](5.2,3.8){$\beta$}
 \psline[linestyle=dashed](5.8,3.6)(5.8,4.5)
 \psline(4.2,3.6)(4.2,4.5)
 \psline(4.6,3.6)(4.6,4.5)
 \psline(5,3.6)(5,4.5)
 \psline(4.2,4.6)(4.4,5.3)
 \psline(4.6,4.6)(4.8,5.3)
 \psline(5  ,4.6)(5.2,5.3)
 \psline(5.2,4.6)(5.4,5.3)
 \psline(5.4,4.6)(5.6,5.3)
 \psline[linestyle=dashed](5.6,4.6)(5.8,5.3)
 \psline[linestyle=dashed](5.8,4.6)(5.6,4.8)(5,4.95)(4.4,5.1)(4.2,5.3)
 \rput[b](6.4,4.7){$\tau_{n+m}$}
 \rput[b](7.3,3){$=$}
 \psline[linestyle=dashed](8.55,2)(8.1,2.9)
 \psline(8.55,2)(8.4,2.9)
 \psline(8.55,2)(8.7,2.9)
 \psline[linestyle=dashed](8.55,2)(9,2.9)
 \rput[b](9.7,2){$\tau_{n+1}(\beta)$}
 \psline[linestyle=dashed](8.1,3)(8.1,3.5)
 \psline[linestyle=dashed](8.1,3.5)(7.8,4)
 \psline(8.1,3.5)(8  ,4)
 \psline(8.1,3.5)(8.2,4)
 \psline(8.1,3.5)(8.4,4)
 \rput[b](8.5,3.5){$\alpha$}
\end{pspicture}
\]
where the last equality follows form equation \eqref{compos_cyclic3}, $\alpha\in \widehat{\mathcal O}(\vec X;\varnothing) =\mathcal O(m)$ has $m+1$ inputs, and $\beta\in\widehat{\mathcal
O}(\vec Y;\da)=\mathcal O(n)$ (or similarly $\beta\in\widehat{\mathcal
O}(\vec Y;\f)$) has $n$ inputs. It is clear that this will satisfy
equivariance, since equivariance was just used to define the
composition. The next lemma establishes the final property for
$\widehat{\mathcal O}$ being a 0/1-operad.
\begin{Lem}
The composition in $\widehat{\mathcal O}$ satisfies the
associativity axiom.
\end{Lem}
\begin{proof} By definition the composition is
just the usual composition in $\mathcal O(n)$, except for inserting
trees in the last input of elements in $\widehat{\mathcal
O}(\vec X;\varnothing)$. Thus, except for composition in the last
spot, associativity of $\widehat{\mathcal O}$ follows from the
associativity of $\mathcal O$.

Now, let $\alpha\in \widehat{\mathcal O}(\vec X;\varnothing)\cong
\mathcal O(m)$, $\beta\in \widehat{\mathcal O}(\vec Y;y)\cong
\mathcal O(n)$, and $\gamma\in \widehat{\mathcal O} (\vec Z;z) \cong
\mathcal O(p)$, where $y,z\in\{\f,\da\}$. Then, associativity is
satisfied, because for $1\leq j\leq m$, it is
\begin{multline*}
(\alpha\circ_{m+1}\beta)\circ_{j} \gamma
\stackrel{\mathit{\eqref{def_cyclic_compos}}}{=} (\tau_{n+1}(\beta)
\circ_1 \alpha)\circ_j \gamma \stackrel{\mathit{op.comp}}{=}\\
=\tau_{n+1}(\beta) \circ_1 (\alpha\circ_j \gamma )
\stackrel{\mathit{\eqref{def_cyclic_compos}}}{=} (\alpha\circ_j
\gamma )\circ_{m+p}\beta,
\end{multline*}
and for $m< j< m+n$, it is
\begin{multline*}
(\alpha\circ_{m+1}\beta)\circ_{j} \gamma
\stackrel{\mathit{\eqref{def_cyclic_compos}}}{=} (\tau_{n+1}(\beta)
\circ_1 \alpha)\circ_j \gamma  \stackrel{\mathit{op.comp}}{=}
(\tau_{n+1}(\beta) \circ_{j-m+1} \gamma)\circ_1 \alpha
\stackrel{\eqref{compos_cyclic2}}{=}\\ = \tau_{n+p}(\beta
\circ_{j-m}\gamma)\circ_1\alpha \stackrel{\mathit{
\eqref{def_cyclic_compos}}}{=} \alpha\circ_{m+1}
(\beta\circ_{j-m}\gamma),
\end{multline*}
while
\begin{multline*}
(\alpha\circ_{m+1}\beta)\circ_{m+n} \gamma
\stackrel{\mathit{\eqref{def_cyclic_compos}}}{=} (\tau_{n+1}(\beta)
\circ_1 \alpha)\circ_{m+n} \gamma  \stackrel{\mathit{
\eqref{def_cyclic_compos}}}{=}\\ = \tau_{p+1}(\gamma)\circ_1
(\tau_{n+1}(\beta) \circ_1 \alpha) \stackrel{\mathit{op.comp}}{=}
(\tau_{p+1}(\gamma)\circ_1 \tau_{n+1}(\beta)) \circ_1 \alpha
\stackrel{\eqref{compos_cyclic3}}{=}\\ = \tau_{n+p}(\beta
\circ_{n}\gamma) \circ_1 \alpha \stackrel{\mathit{
\eqref{def_cyclic_compos}}}{=} \alpha\circ_{m+1} (\beta
\circ_{n}\gamma).
\end{multline*}
\end{proof}
\end{Def}

We end this section by giving a presentation of $\widehat{\mathcal O}$ in terms of generators and relations, when $\mathcal O$ is given by quadratic generators and relations. Let us first recall the notion of operads given by generators and relations.

\begin{Def}
Fix a set of colors $C$. Then let $E=\{E^{x,y}_z\}_{x,y,z\in C}$ be a collection of $k$-vector
spaces, together with an $S_2$-action compatible with the colors.
We want $E$ to be the binary generating set of a colored operad,
where $x$, $y$ and $z$ correspond to the colors of the edges of a binary
vertex, i.e. a vertex with exactly two incoming and one outgoing edge. Let $T$ be a rooted, colored tree where each vertex is binary, and let $v$ be a vertex with colors $(x,y;z)$ in $T$. Then, we define $E(v):= \left( E^{x,y}_z \oplus E^{y,x}_z \right)_{S_2}$,
and with this, we set $E(T):=\bigotimes_{\text{vertex }v\text{ of }T} E(v)$.
We define the free colored operad $\mathcal F(E)$ generated by $E$ to be
\[
\mathcal F(E)(\vec X;z):= \bigoplus_{
  \text{binary trees $T$ of type }(\vec X;z)
} E(T).
\]
The $S_n$-action is given by an obvious permutation of the leaves of the tree using
the $S_2$-action on $E$, and the composition maps are given by
attaching trees. This definition can readily be seen to define a
colored operad.

An ideal $\mathcal I$ of a colored operad $\mathcal P$ is a
collection of $S_n$-sub-modules $\mathcal I(\vec X;z)\subset
\mathcal P(\vec X;z)$ such that $f\circ_{i} g$ belongs to the ideal
whenever $f$ or $g$ or both belong to the ideal. A colored operad $\mathcal P$ is said to be quadratic if $\mathcal P=\mathcal F(E)/(R)$ where $\mathcal F(E)$ is the free
colored operad on some generators $E$, and $(R)$ is the ideal in
$\mathcal F(E)$ generated by a subspace with 3 inputs, called the
relations, $R\subset \bigoplus_{w,x,y,z\in C}
\mathcal F(E)(w,x,y;z)$.
\end{Def}

We recall from \cite[(3.2)]{GeK}, that an operad $\mathcal O$
is called cyclic quadratic if it is quadratic, with generators
$E$ and relations $R$, so that the $S_2$-action on $E$ is naturally
extended to a $S_3$-action via the sign-representation $sgn:S_3\to S_2$, and
$R\subset\mathcal F(E)(3)$ is an $S_4$-invariant subspace. In this
case, $\mathcal O$ becomes a cyclic operad, see \cite[(3.2)]{GeK}.

The following lemma is straight forward to check.

\begin{Lem}\label{O_hat_quadratic}
Let $\mathcal O$ be cyclic quadratic with generators $E$
and relations $R\subset\mathcal F(E)(3)$. Then $\widehat{\mathcal O}$
is generated by $\widehat{E}:=\widehat{E} ^{\ff,\ff}_{\ff}\oplus \widehat{E} ^{\ff,\dda}_{\dda} \oplus \widehat{E} ^{\dda,\ff}_{\dda} \oplus \widehat{E} ^{\dda,\dda}$, defined as
\begin{eqnarray*}
\widehat{E} ^{\ff,\ff}_{\ff}:=E&\subset \widehat {\mathcal
O}(\f,\f;\f),\\%
\widehat{E} ^{\ff,\dda}_{\dda}:=E&\subset \widehat {\mathcal
O}(\f,\da;\da) ,\\%
\widehat{E} ^{\dda,\ff}_{\dda}:= E&\subset \widehat {\mathcal
O}(\da,\f;\da),\\%
\widehat{E}^{\dda,\dda}:=k&\subset \widehat {\mathcal O}(\da,\da;\varnothing),
\end{eqnarray*}
and has relations
\begin{eqnarray*}
R\subset\mathcal F(E)(3)\cong\mathcal F(\widehat{E})(\f,\f,\f;\f), \\
R\subset\mathcal F(E)(3)\cong\mathcal F(\widehat{E})(\f,\f,\da;\da), \\
R\subset\mathcal F(E)(3)\cong\mathcal F(\widehat{E})(\f,\da,\f;\da), \\
R\subset\mathcal F(E)(3)\cong\mathcal F(\widehat{E})(\da,\f,\f;\da).
\end{eqnarray*}
together with the relations
\begin{eqnarray*}
G\subset \mathcal F(\widehat{E})(\da,\da,\f;\varnothing), \\ G\subset
\mathcal F(\widehat{E})(\da,\f,\da;\varnothing), \\ G\subset \mathcal
F(\widehat{E})(\f,\da,\da;\varnothing),
\end{eqnarray*}
where $G$ corresponds for a given coloring to the space
\[ G:= span \left<
\begin{pspicture}(0,0.2)(1,1)
 \psline(0.5,0)(0.3,0.5)
 \psline(0.5,0)(0.7,0.5)
 \psline(0.3,0.5)(0.3,1)
 \psline(0.3,0.7)(0.1,1)
 \psline(0.3,0.7)(0.5,1)
 \rput[b](0.7,0.8){$\alpha$}
\end{pspicture}
-
\begin{pspicture}(0,0.2)(2,1)
 \psline(0.5,0)(0.3,0.5)
 \psline(0.5,0)(0.7,0.5)
 \psline(0.7,0.5)(0.7,1)
 \psline(0.7,0.7)(0.9,1)
 \psline(0.7,0.7)(0.5,1)
 \rput[b](1.5,0.6){$\tau_3(\alpha)$}
\end{pspicture}
\text{ , for all } \alpha\in E^{\ff,\ff}_{\ff} \right>.
\]
\end{Lem}

\section{Koszulness of $\widehat{\mathcal O}$}
\label{quadrat-koszul}

This section is concerned with our main theorem, that Koszulness for $\mathcal O$ implies Koszulness for $\widehat{\mathcal O}$. To set up notation, we briefly recall the notion of quadratic dual, cobar dual and Koszulness of a (colored) operad.

Recall that if the vector space $V$ is an $S_n$-module and $sgn_n$ is the
sign representation, then we defined $V^\vee$ to be $V^*\otimes sgn_n$, where $V^*=Hom(V,k)$ denotes the dual space.
\begin{Def}
Let $C$ be a set of colors. For every quadratic colored operad $\mathcal P$, we define
the quadratic dual colored operad $\mathcal P^!:=\mathcal
F(E)^\vee/(R^\perp)$, where $(R^\perp)$ is the ideal in $\mathcal
F(E)^\vee$ generated by the orthogonal complement $R^\perp$ of $R$
as a subspace of $\left(\bigoplus_{w,x,y,z\in C} \mathcal F(E)(w,x,y;z) \right) ^\vee$. Notice that $\mathcal F(E)(\vec X;x)^\vee = \mathcal F(E^\vee) (\vec X;x)$, so that $\mathcal P^!$ is generated by $E^\vee$ with relations $R^\perp$, see \cite[(2.1.9)]{GK}.

Now if $\mathcal P=\mathcal F(E)/(R)$ is a quadratic colored operad,
then we can follow the definition from \cite[(3.2.12)]{GK}, and define the cobar dual colored operad $\textbf{D} (\mathcal P)$ of $\mathcal P$, to be given by the complexes $\textbf{D} (\mathcal P)(\vec X;x)$ concentrated in non-positive degree, 
$\textbf{D} (\mathcal P)(\vec X;x):=$
\begin{multline*}
\bigoplus_{
\substack{
  \text{trees  $T$ of}\\
  \text{type $(\vec X;x)$},\\
  \text{no internal edge}}}
\mathcal P(T)^*\otimes \Det(T)\stackrel{\de}{\lora}
\bigoplus_{
\substack{
  \text{trees } T \text{ of}\\
  \text{type }(\vec X;x),\\
  \text{1 internal edge}}}
\mathcal P(T)^*\otimes \Det(T)\stackrel{\de}{\lora}\\
\lora \dots\stackrel{\de}{\lora} \bigoplus_{
\substack{
  \text{trees } T \text{ of}\\
  \text{type }(\vec X;x),\\
  \text{binary tree}}}
\mathcal P(T)^*\otimes \Det(T).
\end{multline*}
Here, $\mathcal P(T)^*$ denotes the dual of $\mathcal P(T)$, and $\Det(T)$ denotes the top exterior power on the space $k^{Ed(T)}$, where $Ed(T)$ is the space of edges of the tree $T$. By definition, we let the furthest right space whose sum is over binary trees, be of degree zero, and all other spaces be in negative degree. 

In general, the zero-th homology of this complex is always canonically isomorphic to the quadratic dual $\mathcal P^!$, i.e. $H^0(\textbf{D}(\mathcal P)(\vec X;x))\cong \mathcal P^!(\vec X;x)$, see \cite[(4.1.2)]{GK}. The quadratic operad $\mathcal P$ is then said to be Koszul if the cobar dual on the quadratic dual $\textbf{D}(\mathcal P^!)(\vec X;x)$ is quasi-isomorphic to $\mathcal P(\vec X;x)$, i.e., by the above, $\textbf{D}(\mathcal P^!)(\vec X;x)$ has homology concentrated in degree zero.
\end{Def}

We now state our main theorem.
\begin{Thm}\label{O_hat_Koszul}
If $\mathcal O$ is cyclic quadratic and Koszul, then  $\widehat{\mathcal O}$ has a resolution, which for a given sequence $\vec X$ of colors, with $|\vec X|=n$, is given by the quasi-isomorphisms
\begin{eqnarray*}
\textbf{D}(\widehat{\mathcal O^!}) (\vec X;\f)&\to& \widehat{
\mathcal O}(\vec X;\f)=\mathcal O(n) \\
\textbf{D}(\widehat{\mathcal O^!}) (\vec X;\da)&\to& \widehat{
\mathcal O}(\vec X;\da)=\mathcal O(n) \\
\textbf{D}(\widehat{\mathcal O^!}) (\vec X;\varnothing) &\to&
\widehat{ \mathcal O}(\vec X;\varnothing)=\mathcal O(n-1)
\end{eqnarray*}
\end{Thm}
\begin{proof}
Koszulness of $\mathcal O$ means exactly that the first and second maps are quasi-isomorphisms. The proof that the third map is also a quasi-isomorphism will concern the rest of this section.

We need to show that the homology of $\textbf{D}(\widehat{\mathcal O^!}) (\vec X;
\varnothing)$ is concentrated in degree $0$:
\begin{eqnarray}\label{H0}
H_0 \left(\textbf{D}(\widehat{\mathcal O^!}) (\vec X;\varnothing)\right
)&=&\widehat{ \mathcal O}(\vec X;\varnothing)\\ %
\label{Hi<0} H_{i<0} \left(\textbf{D}(\widehat{\mathcal O^!})
(\vec X;\varnothing)\right)&=&\{0\}
\end{eqnarray}
As mentioned in definition \ref{def_O_hat}, it is enough to restrict attention to the case where $\vec X = (\da,\f,\ldots,\f,\da)$. Let us first show the validity of equation \eqref{H0}. The map $\textbf{D} (\widehat{\mathcal O^!}) (\vec X;\varnothing)\to \widehat{ \mathcal O}(\vec X;\varnothing)$ expanded in degrees of $\textbf{D}(\widehat{\mathcal O^!})(\vec X;\varnothing)$ is written as 
$\cdots \stackrel{\de}{\lora} \textbf{D}(\widehat{\mathcal O^!})
(\vec X;\varnothing)^{-1} \stackrel{\de}{\lora}
\textbf{D}(\widehat{\mathcal O^!}) (\vec X;\varnothing)^0
\stackrel{proj}\lora \widehat{\mathcal O}(\vec X;\varnothing)$.
As $\mathcal O$ and thus $\mathcal O^!$ are quadratic, we have the
following identification, using the language and results of lemma
\ref{O_hat_quadratic}:
\begin{eqnarray*}
\widehat{\mathcal O}(\vec X;\varnothing)&=&\mathcal F(
\widehat{E})/(R,G) (\vec X;\varnothing) \\
\textbf{D}(\widehat{\mathcal O^!}) (\vec X;\varnothing)^0&=&
\bigoplus_{ \substack{
  \text{binary trees }T \\
  \text{of type }(\vec X;\varnothing)
}} (\widehat{E^\vee}(T))\otimes \Det(T) =\mathcal
F(\widehat{E})(\vec X;\varnothing)\\ \textbf{D}(\widehat{\mathcal
O^!}) (\vec X;\varnothing)^{-1}&=&\bigoplus_{ \substack{
  \text{trees }T\text{ of type }(\vec X;\varnothing), \\
  \text{binary vertices, except}\\
  \text{one ternary vertex}
}} \left(\widehat{\mathcal O^!}(T)\right)\otimes \Det(T) \\ &=&
\left\{
\begin{array}{c}
\text{space of relations in $\mathcal F(\widehat{E})(\vec X;\varnothing)$}\\
\text{generated by $R$ and $G$}
\end{array}
\right\}
\end{eqnarray*}
The last equality follows, because the inner product relations for
$\mathcal O^!$, which are the relation space $G$ for the cyclic
quadratic operad $\mathcal O^!$ from lemma \ref{O_hat_quadratic},
are the orthogonal complement of the inner product relations for
$\mathcal O$. Hence, the map $proj$ is surjective with
kernel $\de\big(\textbf{D}(\widehat{\mathcal O^!}) (\vec X;
\varnothing )^{-1}\big)$. This implies equation \eqref{H0}.

As for equation (\ref{Hi<0}), we will use an induction that shows
that every closed element in $\textbf{D}(\widehat{\mathcal O^!})
(\vec X;\varnothing)^{-r}$, for $r\geq 1$ is also exact. The argument 
will use an induction which slides all of the ``full'' inputs from one of
the two ``dashed'' inputs to the other. As a main ingredient of this, we will employ the products $*$ and $\#$ defined below, which are used to uniquely decompose an element in $\textbf{D}(\widehat {\mathcal O^!})(\vec X;\varnothing)$ as a sum of products of $*$ and $\#$.

We need the following definition. Given two
decorated trees $\varphi\in \textbf{D}(\mathcal O^!)(k), \psi \in
\textbf{D}(\mathcal O^!)(l)$, we define new elements
$\varphi * \psi$ and $\varphi \# \psi$ in
$\textbf{D}(\widehat{ \mathcal O^!})(\da,\f,\ldots,\f,\da;\varnothing)$ as
follows. First, for $\varphi * \psi$ take the outputs of
$\varphi$ and $\psi$ and insert them into the unique inner product
decorated by the generator $1\in \widehat{\mathcal O^!}(\da,\da;\varnothing)$:
\[
\begin{pspicture}(-3,0)(6,4)
 \psline[linestyle=dashed](2.5,0.5)(0,3)
 \psline[linestyle=dashed](2.5,0.5)(5,3)
 \rput(2.7,0.3){\tiny $1$}
 \pscircle[linestyle=dotted](4,2.2){1.4}
 \rput(-0.2,1){$\varphi$}
 \psline(1.5,1.5)(2,3)
 \psline(1.5,1.5)(1,3)
 \psline(1.3,2.1)(1.5,3)
 \psline(0.4,2.6)(0.4,3)
 \psline(0.7,2.3)(0.7,3)
 \rput(1.4,1.3){\tiny $\alpha_1$}
 \rput(0.2,2.4){\tiny $\alpha_2$}
 \rput(0.6,2.1){\tiny $\alpha_3$}
 \rput(1,2.3){\tiny $\alpha_4$}
 \pscircle[linestyle=dotted](1,2.2){1.4}
 \rput(5.2,1){$\psi$}
 \psline(3.5,1.5)(3,3)
 \psline(3.5,1.5)(4,3)
 \psline(4.5,2.5)(4.2,3)
 \psline(3.3,2.1)(3.3,3)
 \psline(3.3,2.1)(3.6,3)
 \rput(3.6,1.2){\tiny $\beta_1$}
 \rput(3,2){\tiny $\beta_2$}
 \rput(4.6,2.2){\tiny $\beta_3$}
 \rput(-1.5,2){$\varphi * \psi =$}
\end{pspicture}
\]
As for $\varphi\# \psi$, we assume, that 
$\varphi\in \textbf{D}(\mathcal O^!)(k)$ with $k\geq 2$. Then
identify $\varphi$ with an element in $\textbf{D} (\widehat{
\mathcal O^!})(\da,\f,\ldots,\f,\da;\varnothing)$ by interpreting the
lowest decoration $\alpha_1\in O^!(m)$ of $\varphi$, as an inner
product $\alpha_1\in \widehat{\mathcal O^!}(\da,\f,\ldots,\f,\da;\varnothing)=\mathcal O^!(m)$. $\varphi\# \psi$ is defined by attaching $\psi$ to the dashed 
input labeled by $\alpha_1$:
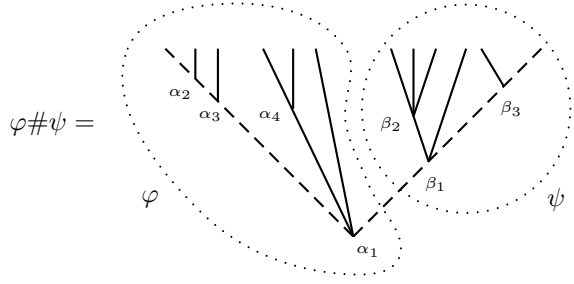
\[
\begin{pspicture}(-3,0)(6,4)
 \psline[linestyle=dashed](2.5,0.5)(0,3)
 \psline[linestyle=dashed](2.5,0.5)(5,3)
 \rput(2.7,0.3){\tiny $\alpha_1$}
 \pscircle[linestyle=dotted](4,2.2){1.4}
 \rput(-0.2,1){$\varphi$}
 \psline(2.5,0.5)(2,3)
 \psline(2.5,0.5)(1.3,3)
 \psline(1.7,2.2)(1.7,3)
 \psline(0.4,2.6)(0.4,3)
 \psline(0.7,2.3)(0.7,3)
 \rput(0.2,2.4){\tiny $\alpha_2$}
 \rput(0.6,2.1){\tiny $\alpha_3$}
 \rput(1.4,2.1){\tiny $\alpha_4$}
 \psccurve[linestyle=dotted](0.5,0.8)(-0.5,3)(1,3.6)(2.4,3)(2.4,2)(2.7,1.2)(3.1,0.2)(2.4,0)
 \rput(5.2,1){$\psi$}
 \psline(3.5,1.5)(3,3)
 \psline(3.5,1.5)(4,3)
 \psline(4.5,2.5)(4.2,3)
 \psline(3.3,2.1)(3.3,3)
 \psline(3.3,2.1)(3.6,3)
 \rput(3.6,1.2){\tiny $\beta_1$}
 \rput(3,2){\tiny $\beta_2$}
 \rput(4.6,2.2){\tiny $\beta_3$}
 \rput(-1.5,2){$\varphi \# \psi =$}
\end{pspicture}
\]
Both $\varphi$ and $\psi$ are elements of $\textbf{D} (\mathcal O^!)$
and thus uncolored. The colorations for $\varphi *\psi$ and $\varphi\#\psi$ are uniquely
determined by having the first and last entry dashed. The operations 
$*$ and $\#$ are extended to $\textbf{D} (\mathcal O^!)\otimes \textbf{D}
(\mathcal O^!)\to \textbf{D} (\widehat{\mathcal O^!})(\vec X; \varnothing)$ by bilinearity.

It is important to notice that every labeled tree, whose first and last inputs are ``dashed", can uniquely be written as a product of $*$ or $\#$. To be more precise, suppose $\varphi$ and $\psi$ are labeled trees with $k+1$ and $l+1$ inputs respectively. We restrict our attention to the case of elements in $\textbf{D} (\widehat{\mathcal O^!})(\vec X; \varnothing)$ whose ``dashed'' inputs are labeled to be the first and the last input, and thus
appear in the planar representation on the far left and the far
right. Let $\sigma$ be a $(k,l)$-shuffle, i.e. a permutation of $\{1,\dots,k+l\}$
such that $\sigma(1)<\dots <\sigma(k)$ and $\sigma(k+1)<\dots <\sigma(k+l)$.
Then define $\phi *_\sigma \psi$, resp. $\phi \#_\sigma \psi$, as the composition
of $*$, resp. $\#$, with $\sigma$ applied to the ``full'' leaves of
the resulting labeled tree. The ``dashed'' inputs remain far left and
far right. With these notations it is now clear, that every labeled tree in $\textbf{D}(\widehat {\mathcal O^!})(\vec X;\varnothing)$ with ``dashed'' first and last inputs, can uniquely be written in the form $\varphi *_\sigma\psi$ or $\varphi\#_\sigma\psi$ for some $\varphi,\psi$ and $\sigma$.

For each $r\geq 1$, we show that the $(-r)$th homology of $\textbf{D}(\widehat{\mathcal O^!})(\vec X;\varnothing)$ vanishes by decomposing every element in $\textbf{D}(\widehat{\mathcal O^!})(\vec X;\varnothing)^{-r}$ as a sum of terms of the form $\varphi *_\sigma\psi$ or $\varphi\#_\sigma\psi$, and then performing and induction on the degree of $\psi$. More precisely, let us define the order of $\varphi*_\sigma \psi$ or $\varphi\#_\sigma \psi$ to be the degree of $\psi$ in $\textbf{D}(\mathcal O^!)$. Then for $s\in \mathbb N$, we claim the following statement:
\begin{itemize}
\item
Let $\chi\in \textbf{D}(\widehat{\mathcal O^!})
(\vec X;\varnothing)^{-r}$ be a closed element $\de(\chi)=0$. Then
$\chi$ is homologous to a sum $\sum_i \sum_\sigma \varphi_i *_\sigma
\psi_i+\sum_j \sum_\sigma \varphi'_j \#_\sigma\psi'_j$, where the order of each term 
is less or equal to $(-s)$.
\end{itemize}
Rather intuitively this means, that for smaller $(-s)$, $\chi$ is homologous to
decorated trees whose total degree is more and more concentrated on
the right branch of the tree.

It is easy to the that the above is true for $s=1$. As for the inductive step,
let $\chi=\sum_i \sum_\sigma \varphi_i *_\sigma \psi_i+\sum_j\sum_\sigma  \varphi'_j
\#_\sigma\psi'_j$, where we assume an expansion so that
$\{\varphi_i\}_i$ are linear independent, $\{\varphi'_j \}_j$ are
linear independent, but the $\psi_i$ and the $\psi'_j$ are allowed to
be linear combinations in $\textbf{D}(\mathcal O^!)$. We claim that
those elements $\psi_i$ and $\psi'_j$, which are of degree $-s$, are
closed in $\textbf{D} (\mathcal O^!)$. This follows from
$\de(\chi)=0$ and the inductive hypothesis, because the only terms of
the boundary $\de(\chi)$, which are of order $-s+1$, are terms of the form
$\varphi_i *_\sigma \de(\psi_i)$ and $\varphi'_j \#_\sigma
\de(\psi'_j)$. But then, $\psi_i$ and $\psi'_j$ are necessarily of degree
$-s$. The exactness of $\textbf{D} (\mathcal O^!)$ at the degree $-s$
implies that $\psi$s of the degree $-s$ are exact, i.e. there are elements $\Psi_i, \Psi'_j\in \textbf{D}(\mathcal O^!)$ such that $\psi_i=\de(\Psi_i)$, $\psi'_j =\de(\Psi'_j)$. Then, take
$S:=\de\left( \sum_i\sum_\sigma \varphi_i *_\sigma \Psi_i+\sum_j
\sum_\sigma \varphi'_j \#_\sigma\Psi'_j\right)$, where the sum is
over those $i$'s and $j$'s which have the constructed $\Psi_i$'s and
$\Psi'_j$'s. As the total degree of the $\Psi$'s is $-s-1$, we see
that the only terms of $S$ with degree greater or equal to $-s$ are
the terms
$$\sum_i \sum_\sigma \varphi_i *_\sigma \de(\Psi_i) +
\sum_j\sum_\sigma \varphi'_j \#_\sigma\de(\Psi'_j)
=\sum_i \sum_\sigma \varphi_i *_\sigma \psi_i+
\sum_j \sum_\sigma \varphi'_j \#_\sigma\psi'_j.
$$
It follows that $\chi-S$ only contains terms of degree less or equal
to $-s-1$, and $\chi$ is homologous to $\chi-S$. This concludes the inductive step.

We complete the proof of the theorem by noticing that $\textbf{D}(\mathcal O^!)
(l)$ is concentrated in finite degrees, so that the $\psi_i$s and
$\psi'_j$s eventually have to be identically $0$. As a consequence every closed
element $\chi$ is eventually homologous to $0$.
It follows that the complex $ \textbf{D}(\widehat{\mathcal O^!})
(\vec X;\varnothing)$ has no homology in degrees $r\neq 0$.
\end{proof}

\section{Algebras over $\widehat{\mathcal O}_\infty$}
\label{homotop-ip}

In this section, we investigate algebras over $\widehat{\mathcal O}_\infty:=\textbf{D} (\widehat{\mathcal O^!})$. The particular cases of the associative operad $\mathcal Assoc$ and the commutative operad $\mathcal Comm$ will be considered.

Before looking at $\widehat{\mathcal O}_\infty$, let us first consider algebras over $\widehat{\mathcal O}$. These are given by ``algebra maps'' $\mathcal O(n) \otimes A^{\otimes n}\to A$, ``module maps'' $\bigoplus_{r+s=n-1} \mathcal O(n) \otimes A^{\otimes r}\otimes M\otimes A^{\otimes s}\to M$, and since $\widehat{\mathcal O}(\vec X;\varnothing) =\mathcal O(n-1)$ for $|\vec X|=n$, we also have ``inner product maps''
$ \mathcal O(n-1) \to Hom(A^{\otimes i-1} \otimes M
\otimes A ^{\otimes j-i-1} \otimes M \otimes A^{\otimes n-j} ,k)$.
Notice that in the lowest case $n=2$, the ``inner product map" $\mathcal O(1)\to Hom(M\otimes M,k)$ determines a map $<.,.>:M\otimes M\to k$, given by the image of the unit $1_k\in k=\mathcal O(1)$. Note that $<.,.>$ is invariant under the module maps mentioned above, and using the composition and the $S_n$-action of $\widehat{\mathcal O}$,  all the higher inner product maps are determined by $<.,.>$ together with the module maps.

We now describe algebras over the operad $\textbf{D} (\widehat{\mathcal O^!})$. This will be given in terms of derivations and maps of free $\mathcal O$-algebras and free $\mathcal O$-modules. 
\begin{Def}
Recall from \cite[(1.3.4)]{GK}, that a free $\mathcal P$-algebra is generated by the $k$-vector space $A$ given by $F_{\mathcal P} A:=\bigoplus_{n\geq 1}(\mathcal P(n)\otimes A^{\otimes n})_{S_n}$. $F_{\mathcal P} A$ is an algebra over $\mathcal P$, i.e., there are maps $ \gamma:\mathcal P(n)\otimes (F_{\mathcal P} A )^{\otimes n} \to F_{\mathcal P} A$ coming from the composition in $\mathcal P$ and the tensor products, which satisfy the required
compatibility conditions, see \cite[(1.3.2)]{GK}. An algebra derivation $d\in \mathrm{Der}(F_{\mathcal P}A)$ is defined to be a map from $F_{\mathcal P} A$ to itself, making the following
diagram commute:
\[
\begin{CD}
\mathcal P(k) \otimes (F_{\mathcal P}A)^{\otimes k}& @>\gamma>> &
F_{\mathcal P}A\\ @V\sum_i \mathrm{id}\otimes \mathrm{id}^{\otimes i}
\otimes d \otimes \mathrm{id}^{\otimes (k-i-1)}VV & & @VVdV\\
\mathcal P(k) \otimes (F_{\mathcal P}A)^{\otimes k}& @>\gamma>> &
F_{\mathcal P}A
\end{CD}
\]

Similarly, if $A$ and $M$ are $k$-vector spaces, we define the free module $M$ over $A$ to be 
$F_{\mathcal P,A} M := \bigoplus_{n\geq 1}
\left(\bigoplus_{r+s=n-1} \mathcal P(n) \otimes A^{\otimes r}
\otimes M \otimes A^{\otimes s}\right)_{S_{n}}$.
Then $F_{\mathcal P,A} M$ is a module over $F_{\mathcal P}A$ over
$\mathcal P$, which means that there are maps $\gamma^M:\bigoplus_{r+s=n-1}
\mathcal P(n)\otimes (F_{\mathcal P} A )^{\otimes r} \otimes
F_{\mathcal P,A} M \otimes (F_{\mathcal P} A )^{\otimes s} \to
F_{\mathcal P,A} M$ given by composition of elements of the operad
and tensor product of elements of $A$. These maps satisfy the
required module axioms, see \cite[(1.6.1)]{GK}. A module derivation $g\in \mathrm{Der}_d (F_{\mathcal P, A}M)$ over $d\in \mathrm{Der}(F_{\mathcal P}A)$, is define to be a map from $F_{\mathcal P, A}M$ to itself, making the following diagram commutative:
\[
\begin{CD}
\mathcal P(k) \otimes (F_{\mathcal P}A)^{\otimes (k-1)}\otimes
F_{\mathcal P, A}M & @>\gamma^M>> & F_{\mathcal P, A}M\\ @V\sum_{i<k} \mathrm{id}\otimes \mathrm{id}^{\otimes i}
\otimes d \otimes \mathrm{id}^{\otimes (k-i-1)}+ \mathrm{id}\otimes
\mathrm{id}^{\otimes (k-1)} \otimes g VV
& & @VVgV\\ \mathcal P(k) \otimes (F_{\mathcal P}A)^{\otimes
(k-1)}\otimes F_{\mathcal P, A}M & @>\gamma^M>> & F_{\mathcal P, A}M
\end{CD}
\]

Finally, a module map $f\in \mathrm{Mod}(F_{\mathcal P, A}M, F_{\mathcal
P,A}N)$ is defined to be a map from $F_{\mathcal P,A}M$ to $F_{\mathcal P,A}N$
making the following diagram commutative:
\[
\begin{CD}
\mathcal P(k) \otimes (F_{\mathcal P}A)^{\otimes (k-1)}\otimes
F_{\mathcal P, A}M & @>\gamma^M>> & F_{\mathcal P, A}M\\
@V\mathrm{id}\otimes \mathrm{id}^{\otimes (k-1)} \otimes f VV & &
@VVfV\\ \mathcal P(k) \otimes (F_{\mathcal P}A)^{\otimes
(k-1)}\otimes F_{\mathcal P, A}N & @>\gamma^N>> & F_{\mathcal P, A}N
\end{CD}
\]
\end{Def}

If $\mathcal P$ is a cyclic operad, then we can use this extra datum to associate to every derivation $g$ over a free module $M$ a derivation $h$ over the free module on the dual space $M^\ast=Hom(M,k)$.
\begin{Def}\label{dual-module}
Let $\mathcal P$ be a cyclic operad, and let $A$ and $M$ be a vector space over $k$, which are
finite dimensional in every degree. Assume furthermore, that there are derivations  $d\in \mathrm{Der} (F_{\mathcal P}\,A)$, and $g\in \mathrm{Der}_d (F_{\mathcal P,A}M)$ over $d$. The maps $d$ and $g$ are determined by maps
\begin{eqnarray*}
d_n:\,A&\to& \bigoplus_n \mathcal P (n)\otimes A^{\otimes n},\\ %
g_n:\,M&\to& \bigoplus_{k+l=n-2} \mathcal P (k+l+1)\otimes
A^{\otimes k} \otimes \,M \otimes A^{\otimes l}.
\end{eqnarray*}
Then there is an induced derivation on $M^*$ over $d$ in the
following way. Define $h\in \mathrm{Der}_d (F_{\mathcal P,A}M^*)$ by maps
$h_n:\,M^* \to \bigoplus_{k+l=n-2} \mathcal P(k+l+1)\otimes
A^{\otimes k} \otimes M^* \otimes A^{\otimes l}
$, which are given by its application to $a^\ast_1,\dots,a^\ast_{k+l} \in \,A^\ast$, $m^*_1\in \,M^*$, $m_2\in \,M$ as
\begin{multline*}
h_n(m^*_1) (a^*_1,\dots,a^*_{k},m_2,a^*_{k+1},\dots,a^*_{k+l}) \\ :=(-1)^\epsilon
\cdot \tau_{k+l+2}^{l+1}(g_n(m_2)(a^*_{k+1},\dots,a^*_{k+l},m^*_1,a^*_1,\dots,
a^*_{k})) \in\mathcal P(k+l+1).
\end{multline*}
\[
\begin{pspicture}(0,-1)(5,2.5)
 \psline(2,1)(0.8,0)
 \psline(2,1)(1.2,0)
 \psline(2,1)(1.6,0)
 \psline(2,1)(2.4,0)
 \psline(2,1)(2.8,0)
 \psline(2,1)(3.2,0)
 \psline[arrowsize=0.2, arrowinset=0, linestyle=dashed](2,2)(2,1)
 \rput[l](2.5,1.35){$h$} \rput[m](2,2.2){$m^*_1$}
 \psline[linestyle=dashed](2,0)(2,1)
 \rput(2,-0.4){$a^*_{k+l}\dots a^*_{k+1}$ $m_2$ $a^*_{k} \dots a^*_{1}$\quad}
 \psline[linearc=0.3, arrowsize=0.15]{->}(1.6,1.2)(1.9,1.5)(2.4,1)(1.9,0.5)(1.5,0.9)
 \rput(4.5,1.5){$:=(-1)^\epsilon \cdot$}
\end{pspicture}
\begin{pspicture}(0,-1)(4,2.5)
 \psline(2,1)(0.8,0)
 \psline(2,1)(1.2,0)
 \psline(2,1)(1.6,0)
 \psline(2,1)(2.4,0)
 \psline(2,1)(2.8,0)
 \psline(2,1)(3.2,0)
 \psline[arrowsize=0.2, arrowinset=0, linestyle=dashed](2,2)(2,1)
 \rput[l](2.5,1.35){$\left(\tau_{k+l+2}^{l+1}\otimes
                      id\right) (g)$} \rput[m](2,2.2){$m_2$}
 \psline[linestyle=dashed](2,0)(2,1)
 \rput(2,-0.4){\quad $a^*_{k} \dots a^*_{1}$  $m^*_1$ $a^*_{k+l} \dots a^*_{k+1}$}
 \psline[linearc=0.3, arrowsize=0.15]{->}(1.6,1.2)(1.9,1.5)(2.4,1)(1.9,0.5)(1.5,0.9)
\end{pspicture}
\]
Here $\tau_{k+l+2}^{l+1}$ denotes the $(l+1)$st iteration of
$\tau_{k+l+2}$, and $\epsilon=(|m^*_1|+|a^*_1|+\dots+|a^*_k|)\cdot
(|m_2|+|a^*_{k+1}|+\dots+|a^*_{k+l}|)$. 

A straightforward check shows that if $d^2=0$ and $g^2=0$, then it is also $h^2=0$.
\end{Def}

The following proposition identifies algebras over $\textbf{D} (\widehat{\mathcal O^!})$. Its proof is similar to that in \cite[Proposition (4.2.14)]{GK}.  
\begin{Prop}\label{O_hat_algebras} Let $\mathcal O$ be a cyclic quadratic operad, and let $A$ and $M$ be graded vector spaces over $k$, which are finite dimensional in every degree. Then giving $(A,M)$ the structure of an algebra over $\textbf{D} (\widehat{\mathcal O^!})$ is equivalent to the following data:
\begin{itemize}
\item
a derivation $d\in \mathrm{Der}(F _{\mathcal O^!}\,A^*[1])$ of degree
$1$, with $d^2=0$,
\item
a derivation $g\in \mathrm{Der}_d (F_{\mathcal O^!,\,A^*[1]}M^*[1])$ over
$d$ of degree $1$, with $g^2=0$, which by definition \ref{dual-module} also implies a derivation $h\in \mathrm{Der}_d (F_{\mathcal O^!,\,A^*[1]}M[1])$ over $d$ with $h^2=0$,
\item
a module map $f\in \mathrm{Mod}(F_{\mathcal O^!, A^*[1]}M[1], F_{\mathcal
O^!,A^*[1]}M^*[1])$ of degree $0$ such that $f\circ h = g \circ f$, and satisfying the following symmetry condition: let $f$ be given by maps $f_n:\,M[1]\to \bigoplus_{k+l=n-2} \mathcal O^!(k+l+1)\otimes A^*[1]^{\otimes k}\otimes \,M^*[1] \otimes A^*[1]^{\otimes l}$, then
\begin{multline}
\label{eq:symm}
 f_n(m_2)(\alpha^*;a_1,\dots,a_i,m_1,a'_1,\dots,a'_j)=\\
= (-1)^\epsilon f_n(m_1)(\alpha^*\circ \tau_{i+j+2}^{j+1}
;a'_1,\dots,a'_j,m_2,a_1,\dots,a_i).
\end{multline}
\end{itemize}
where $\epsilon=(|m_2|+|a_1|+\dots+|a_i|+i+1)\cdot (|m_1|+ |a'_1|+\dots+
|a'_j|+j+1)$.
\end{Prop}

If $\mathcal O$ is cyclic quadratic and Koszul, then, by theorem
\ref{O_hat_Koszul}, $\textbf{D}( \widehat{\mathcal O^!} )\cong
\widehat{\mathcal O}$ and we call $(A,M)$ a homotopy $\mathcal O$-algebra with homotopy $\mathcal O$-module and homotopy $\mathcal O$-inner product
if there are derivations $d$ and $g$ together with a module map $f$ satisfying the conditions from proposition \ref{O_hat_algebras}.

\begin{Exa}\label{assoc}
Let $\mathcal Assoc$ be the associative operad, see \cite[(1.3.7)]{GK} and \cite[(2.3)]{GeK}. In this case, $\mathcal Assoc^!=\mathcal Assoc$, and therefore the free $\mathcal Assoc^!$-algebra is given by the tensor algebra on the underlying vector space, i.e. $F_{\mathcal Assoc^!} A^*[1]=T(A^*[1])=\bigoplus_{n} A^*[1]^{\otimes n}$. Similarly, $F_{\mathcal Assoc^!, A^*[1]}M^*[1]=T^{A^*[1]}{M^*[1]}:=\bigoplus_{k,l} M^*[1]^{\otimes k}\otimes A^*[1]\otimes M^*[1]^{\otimes l}$. We see after dualizing, that $(A,M)$ is a homotopy $\mathcal Assoc$-algebra with homotopy $\mathcal Assoc$-module and homotopy $\mathcal Assoc$-inner product if there are coderivations $D:TA[1]\to TA[1]$ and $G:T^{A[1]} M[1]\to T^{A[1]} M[1]$ and a comodule map $F:T^{A[1]} M[1]\to T^{A[1]} M^*[1]$ satisfying $D^2=0, G^2=0, F\circ H=G\circ F$ and the symmetry condition \eqref{eq:symm}, where $H$ is the induced coderivation on the dual space $M^*[1]$.

Thus, we recover exactly the concept of $\infty$-inner
product over an $A_\infty$-algebra as defined in \cite{Tr}, which
additionally satisfies the symmetry condition \eqref{eq:symm}, coming
from switching the factors of $M$. Furthermore, one can see from
\cite[Lemma 2.14]{Tr2}, that this additional symmetry \eqref{eq:symm}
implies that the $\infty$-inner product is symmetric in the sense of
\cite[Definition 2.13]{Tr2}. In \cite{Tr2}, it was shown that if a symmetric
$\infty$-inner product is also non-degenerate, then it induces a BV-structure on the
Hochschild-cohomology of the given $A_\infty$-algebra. It would be
interesting to have a generalization of this result to any cyclic
operad $\mathcal O$, which amounts to look for a similar ``$\mathcal
O$-BV-structure'' on the homology of the chain complex $
C_n^{\mathcal O}(A)=\left(\mathcal O^!(n)^{\vee}\otimes A^{\otimes n}
\right)_{S_n} $ of an $\mathcal O$-algebra $A$ from \cite[(4.2.1)]{GK}.
\end{Exa}

\begin{Exa}\label{exa-comm}
Let $\mathcal Comm$ be the commutative operad, see \cite[(1.3.8)]{GK} and \cite[(3.9)]{GeK}. It is well known, that $\mathcal Comm^!=\mathcal Lie$ is the Lie operad. Thus, we get $F_{\mathcal Comm^!} A^*[1]=L(A^*[1]):=\bigoplus_n \left(\mathcal Lie(n)\otimes A^*[1]^{\otimes n} \right)_{S_n}$ is the free Lie algebra generated by $A^*[1]$, and $F_{\mathcal Comm^!, A^*[1]}M^*[1]= \bigoplus_n \left(\bigoplus_{k+l=n-1} \mathcal Lie(n) \otimes A^*[1]^{\otimes k} \otimes M^*[1] \otimes A^*[1]^{\otimes l}\right)_{S_{n}}$. It is worth noting, that in \cite[Section 2]{Ginot}, these spaces were already considered as $C_\infty$-algebras with $C_\infty$-modules. We thus have the concept of homotopy $\mathcal Comm$-inner products as module maps between $C_\infty$-modules $M$ and $M^*$.
\end{Exa}


\end{document}